\documentclass[11pt]{article}

\usepackage{amsfonts}
\usepackage{amssymb}
\usepackage{amsmath}
\usepackage{amsthm}
\usepackage{mathrsfs}

\usepackage[final]{showkeys}

\usepackage[subnum]{cases} 

\newtheorem{theorem}{Theorem}
\newtheorem{lemma}[theorem]{Lemma}

\newtheorem{corollary}[theorem]{Corollary}

\theoremstyle{definition}
\newtheorem{remark}[theorem]{Remark}

\newtheorem{hypothesis}{Hypothesis}

\def\eqref#1{(\ref{#1})}

\def\proof{{\parindent0pt {\bf Proof.\ }}}
\def\endproof{\hfill{\hbox{$\square$}}}

\def\weak{\rightharpoonup}
\def\dans{\longrightarrow}

\def\imply{\Longrightarrow}

\let\phi=\varphi
\let\epsilon=\varepsilon

\let\lemref=\ref
\let\Bbb=\mathbb

\begin{document}

\title{%
	A global approach to the \\
	Schr\"odinger-Poisson system: \\
	An existence result in the case\\
	of infinitely many states
}
\author{Otared Kavian\\
\small D\'epartement de Math\'ematiques\\
\small Universit\'e de Versailles\\
\small 45, avenue des \'Etats-Unis\\
\small 78035 Versailles Cedex, France.\\
\small {\tt kavian@math.uvsq.fr} 
\and St\'ephane Mischler\\
\small Ceremade \\
\small Universit\'e de Paris--Dauphine \\
\small place du Mar\'echal de Lattre de Tassigny \\
\small 75116 Paris, France \\
\small {\tt  stephane.mischler@dauphine.fr}}

\maketitle

\begin{abstract}
\noindent In this paper we prove the existence of a solution to a nonlinear Schr\"odinger--Poisson eigenvalue problem in dimension $N$, with $N \leq 6$. Our proof is based on a global approach to the determination of eigenvalues and eigenfunctions which allows us to characterize the complete sequence of eigenvalues and eigenfunctions at once, via a variational approach, and thus  differs from the usual
and less general proofs developed  for similar problems in the literature. 
Our method seems to be new for the determination of the spectrum and eigenfunctions for compact and self-adjoint operators, even in a finite dimensional setting.
\end{abstract}

\noindent{\bf Keywords: }{Schr\"odinger--Poisson, variational methods, unitary operators, Palais--Smale condition.\\

}

\section{Introduction}\label{sec:intro}

In this paper we are interested in the study of the following stationary Schr\"odinger system of equations: for a given bounded Lipschitz domain $\Omega \subset {\Bbb R}^N$, find an infinite sequence $(\lambda_{m}, u_m)_{m\geq1}$ and a potential $V$ satisfying 
\begin{numcases}{}
-\Delta  u_m + {\widetilde V}_{0}u_{m} + V u_m  = \lambda_{m} u_m \quad\mbox{in }\, \Omega  
 \label{eq:1} \\
- \Delta V = \sum_{m=1}^{\infty}\rho_{m}| u_m|^2 \quad\mbox{in }\, \Omega  \label{eq:2}
\end{numcases}
where, in addition, we require that 
\begin{numcases}{}
u_m \in H^1_{0}(\Omega),\; \mbox{ for all } m \geq 1, \quad V \in H^1_{0}(\Omega), \label{eq:BC} \\
(u_m)_{m\geq1} \; \mbox{ is a Hilbert basis for }\, L^2(\Omega). \label{eq:ortho} 
\end{numcases}
Here we assume that the external potential ${\widetilde V}_{0}$ and the positive numbers $(\rho_{m})_{m\geq 1}$ are given. 

The  system \eqref{eq:1}-\eqref{eq:2} appears in the modeling of nanoscale semiconductor devices as part of the so-called ``quantum-kinetic subband model", which itself is a simplified model of the full evolution 3D-Schr\"odinger-Poisson equation, see N. Ben Abdallah \& F. M\'ehats \cite{BenAMehats2004,Mehats2006}. 
In order to improve the cost of numerical simulation of the evolution 3D-Schr\"odinger-Poisson equation and taking advantage of the extreme confinement of the electrons in one direction transverse to the transport directions, one can perform a block diagonalisation of the electron Hamiltonian, thanks to a separation
of the confinement and transport directions. This reduction process leads to replace the 3D-Schr\"odinger-Poisson equation by a system of 1D stationary Schr\"odinger equations (for the confinement) coupled to a 2D equation (for the transport). In that reduced model, the system of 1D stationary Schr\"odinger equations
is nothing but \eqref{eq:1}-\eqref{eq:2} where the parameters $(\rho_{m})_{m\geq 1}$ are the sequence of occupation numbers, which may  depend on time and space, and which are given by the above mentioned 2D equation (for the transport).

On the other hand, as it is outlined in  P. Zweifel \cite{Zweifel}, the Schr\"odinger-Poisson system of equations derives from a quantum transport equation, the Poisson-Wigner system. 
After performing the Wigner transform to the former system one ends up with the following system of evolution equations
\begin{numcases}{}
{\rm i} \, \partial_t \psi_m =  -\Delta  \psi_m + {\widetilde V}_{0}\psi_{m} + V \psi_m    
   \\
   \psi_m(0,x) = \psi_{0m}(x)
   \\
- \Delta V = \sum_{m=1}^{\infty}\rho_{m}| \psi_m|^2,
\end{numcases}
with the condition $(\psi_{0m}|\psi_{0k})_{L^2} = \delta_{km}$ for any $k,m \ge 1$. Seeking standing-wave solutions of the form 
$$
\psi_m(t,x) := {\rm e}^{-{\rm i} \, \lambda_m \, t} \, u_m(x)
$$
leads to the equations \eqref{eq:1}-\eqref{eq:2}. 

Existence of solutions to the  system \eqref{eq:1}-\eqref{eq:2} has been proved in R. Illner, O. Kavian, H. Lange \cite{ILK:1998} in the case when $\rho_m = 0$ 
for any $m \ge 2$, and the same method extends to the case  when $\rho_m = 0$ 
for any $m \ge M$, for a given $M \ge 2$. On the other hand, a different but similar eigenvalue problem has been considered by F. Nier \cite{NierF1990,Nier1993CPDE,Nier1993PRSE}. 

In order to solve this system of equations, taking into account the fact that the family $( u_m)_{m\geq1}$ must be contained in $H^1_{0}(\Omega)$ and, at the same time, has to be a Hilbert basis of $L^2(\Omega)$, we observe the following. Let us consider a fixed Hilbert basis of $L^2(\Omega)$, denoted by $(e_{m})_{m\geq1}$, such that $e_{m} \in H^1_{0}(\Omega)$ for all $m \geq1$. For instance such a basis may be given by the eigenfunctions of the Laplace operator on $H^1_{0}(\Omega)$, that is a family satisfying
\begin{equation}
-\Delta e_{m} = \mu_{m}e_{m},\quad e_{m} \in H^1_{0}(\Omega),\quad 
\int_{\Omega} e_{\ell}(x)\, e_m(x)\,dx = \delta_{\ell m}\, ,
\end{equation}
where the sequence of eigenvalues of the Laplace operator, with Dirichlet boundary conditions, is denoted by $(\mu_{m})_{m\geq 1}$.
Now, saying that the family $(u_m)_{m\geq1}$ satisfies condition \eqref{eq:ortho} means that the linear operator $U$ acting on  $L^2(\Omega)$ and defined by
\begin{equation}
\mbox{for all }\; m\geq1,\qquad Ue_{m} :=  u_m,\qquad Uf := \sum_{m \geq 1} (f|e_{m})  u_m\quad \mbox{for }\, f\in L^2(\Omega),
\end{equation}
is a unitary operator, that is $U^*U = UU^* = I$. Therefore determining the whole family $(u_m)_{m\geq1}$ satisfying equations \eqref{eq:1}--\eqref{eq:ortho} is equivalent to find a linear operator $U$ defined on $L^2(\Omega)$ verifying
\begin{equation}\label{eq:Unitary}
U^*U=UU^*= I, \qquad Ue_{j} \in H^1_{0}(\Omega)\quad\mbox{for }\, j\geq 1,
\end{equation} 
and such that the family $u_{j} := Ue_{j}$ is the family of normalized eigenfunctions  of \eqref{eq:1} where $V$ is given by \eqref{eq:2}.
\medskip

In this paper we give a variational formulation of the system \eqref{eq:1}--\eqref{eq:ortho}, yielding a solution in terms of critical points of a real valued functional defined on a subset of the group of unitary operators. More precisely we define a subset ${\Bbb S}$ of unitary operators on $L^2(\Omega)$ as
\begin{equation}\label{eq:DefS}
{\Bbb S} := \big\{ U : L^2(\Omega) \dans L^2(\Omega) \; ; \; U\; \mbox{satisfies \eqref{eq:Unitary}} \big\},
\end{equation}
and then we define a functional $J$ on ${\Bbb S}$ by setting 
\begin{equation}\label{eq:DefJ0}
J_{0}(U) := \sum_{m \geq1}\rho_{m} \int_{\Omega}\left(|\nabla u_{m}(x)|^2 + {\widetilde V}_{0}(x)u_{m}^2(x)\right)dx,
\end{equation}
and then 
\begin{equation}\label{eq:DefJ}
J(U) := J_{0}(U) + {1\over 2} \int_{\Omega}|\nabla V[U](x)|^2dx,
\end{equation}
where $u_{m}:=Ue_{m}$ and $V[U]$ is the solution of \eqref{eq:2}. The main purpose of this paper is to  show that critical points of $J$ on ${\Bbb S}$ yield solutions of the Schr\"odinger--Poisson system and that the minimum of $J$ is achieved on ${\Bbb S}$.
\medskip

It is clear that in order to define the potential $V[U]$ and the functional $J$ on the mani\-fold ${\Bbb S}$ some conditions must be imposed on the sequence $(\rho_{m})_{m\geq1}$. 

Our main result concerning the system of equations \eqref{eq:1}--\eqref{eq:ortho} is the following (regarding the condition $1 \leq N \leq 6$ see Remark \ref{rem:CondN} below):
\begin{theorem}\label{lem:MainThm}
Let $\Omega \subset {\Bbb R}^N$ be a bounded domain and $1 \leq N \leq 6$. Assume that the sequence $(\rho_{m})_{m \geq 1}$ satisfies
\begin{equation}\label{eq:RhojGrowth}
\rho_{m} >0, \qquad \sum_{m \geq 1}  m^{2/N}\,\rho_{m} < \infty, 
\end{equation}
and 
$$ {\widetilde V}_{0}^+ \in L^1(\Omega), \quad {\widetilde V}_{0}^- \in L^{p_{0}}(\Omega), $$ 
for some $p_{0} > N/2$ (or $p_{0}=1$ if $N=1$). Then the Schr\"odinger--Poisson system of equations \eqref{eq:1}--\eqref{eq:ortho} has a solution obtained as the minimum over ${\Bbb S}$ of the functional $J$ defined in \eqref{eq:DefJ}.
\end{theorem}
\bigskip

In order to give a clear exposition of our global approach to the determination of a system of eigenvectors in terms of unitary operators $U$, in this introduction we give an outline of our approach, getting rid of technicalities inherent to an infinite dimensional Hilbert space, and to nonlinear problems. 

Thus, in a first step, assume that  we are given a finite dimensional (complex) Hilbert space $H$ of dimension $n \geq 2$, its scalar product being denoted by $(\cdot|\cdot)$. If $A : H \dans H$ is a self-adjoint, nonnegative operator (matrix), our aim is to define a procedure in which all the eigenvectors of $A$ are determined at once, to compare with a step by step construction of eigenvalues (and eigenvectors) through the construction of critical values of the Rayleigh quotient
$$ {(Au|u)\over (u|u)}\, , \qquad u\neq 0, $$
by a min-max procedure. To this end, for $u_{1},\ldots,u_{n} \in H$ let us define the functional
$$ F(u_{1},\ldots,u_{n}) := \sum_{j = 1}^n \rho_{j}(Au_{j}|u_{j}), $$
where, as above, we assume that the coefficients $\rho_{j}$ verify
$\rho_{j} >0$, and also, for the sake of simplicity of exposition (see below for the general case), here assume moreover that 
$$ \rho_{i}\neq \rho_{j} \quad\mbox{for }\,  i\neq j. $$
Define the subset (or manifold) $S\subset H^n$ by
$$ S := \big\{\, (u_{1},\ldots,u_{n}) \in H^n \; ; \; (u_{i}|u_{j})=\delta_{ij} \, \big\} $$
We claim that upon maximizing or minimizing  $F$ on the manifold $S \subset H^n$, all the eigenvectors of $A$ can be determined (as a matter of fact, any critical point of $F$ yields such a result). 

\medskip

We begin by observing that an orthonormal basis  $e_{1},\ldots,e_{n}$ of $H$ being given once and for all, the manifold $S$ can be identified with the set of unitary operators  $U$ on $H$ such that $U^*U = UU^* = I$, where $I$ is the identity operator on $H$: indeed it is enough to see $u_{j}$ as the $j$-th column of the matrix representation of $U$, that is to set $u_{j} := Ue_{j}$. Then we see that
$(Au_{j}|u_{j}) = (U^*AUe_{j}|e_{j})$,
and denoting by $D$ the diagonal matrix $D:={\rm diag}(\rho_{1},\ldots,\rho_{n})$, that is the matrix defined by  $De_{j}=\rho_{j}e_{j}$, we check easily that  
$$ F(u_{1},\ldots,u_{n}) = \sum_{j\geq 1} \rho_{j}(Au_{j}|u_{j}) = \sum_{j\geq 1} \rho_{j}(U^*AUe_{j}|e_{j}) = {\rm tr}(DU^*AU), $$
where ${\rm tr}(B)$ denotes the trace of the operator (or matrix) $B$.
Finally, considering for instance the minimization of $F$, this can be reformulated in the following way: 
$$ \mbox{minimize } J(U) := {\rm tr}(DU^*AU)\,\mbox{under the constraint }\, U^*U = I. $$
Clearly $J$ is $C^\infty$ (in fact analytic) and positive on the set
$${\bf U}(n) := {\bf U}(n,H) := \{ U : H \dans H \; ; \; U^*U = I\},$$
which is a smooth and compact manifold: therefore  $J$ achieves its minimum at some  point $U_{0}\in {\bf U}(n)$. Now we have to show that the vectors $u_{j} := U_{0}e_{j}$ are indeed the eigenvectors of $A$.

\medskip

Let $M : H \dans H$ be skew-adjoint (that is $M^* = -M$) and consider the one parameter group $U(t) = \exp(tM)U_{0}$ for $t \in {\Bbb R}$; note that since $M^* = -M$, one has  $\exp(tM)^* = \exp(-tM)$ and thus one checks easily that having $U_{0}^*U_{0}=I$, then for all $t\in{\Bbb R}$ one has  $U(t) \in {\bf U}(n)$, and consequently $J(U_{0}) \leq J(U(t))$ for all $t\in {\Bbb R}$. Now 
$$ J(U(t)) = {\rm tr}(DU_{0}^*\exp(-tM)A\exp(tM)U_{0}) \quad\mbox{and}\quad {d J(U(t))\over dt }_{| t= 0} = 0, $$
so that, after a straightforward calculation, we obtain
$$\begin{cases}
\mbox{for all }\, M\,\mbox{such that }\, M^* = -M, \,\mbox{ we have }\\
{\rm tr}(DU_{0}^*MAU_{0}) = {\rm tr}(DU_{0}^*AMU_{0}). 
\end{cases}$$
Setting
$$ B := U_{0}DU_{0}^* \, , $$
and using the fact that for two given matrices $K,L$ we have ${\rm tr}(KL) = {\rm tr}(LK)$, we observe that ${\rm tr}(DU_{0}^*MAU_{0}) = {\rm tr}(MAB)$, and that
$$ {\rm tr}(DU_{0}^*AMU_{0}) = {\rm tr}(BAM) = {\rm tr}(MBA).$$
Summing up, we conclude that
$$ \mbox{for all }\, M\,\mbox{such that }\, M^* = -M,\quad \mbox{we have }\quad
 {\rm tr}(M(AB-BA)) = 0. $$
Taking $M:= (AB-BA)^*$, we conclude that $BA =AB$, that is
$$ U_{0}DU_{0}^*A = AU_{0}DU_{0}^* . $$
Applying this equality to the vector  $u_{j} := U_{0}e_{j}$, and taking into account the definition of the diagonal operator $D$, we obtain (recall that $U_{0}^*U_{0}=I$)
$$ (U_{0}DU_{0}^*)Au_{j} = AU_{0}DU_{0}^*u_{j} = AU_{0}De_{j}= \rho_{j}AU_{0}e_{j} = \rho_{j}Au_{j}, $$
that is $(U_{0}DU_{0}^*)Au_{j} = \rho_{j}Au_{j}$, which means that $Au_{j}$ is an eigenvector of $U_{0}DU_{0}^*$. This implies that  
$D(U_{0}^*Au_{j}) = \rho_{j}(U_{0}^*Au_{j})$,
and we see that $U_{0}^*Au_{j}$ is an eigenvector of $D$ for the eigenvalue  $\rho_{j}$, which is a simple eigenvalue of $D$, corresponding to the eigenvector $e_{j}$. This means that there exists  $\lambda_{j} \in {\Bbb C}$ such that  $U_{0}^*Au_{j} = \lambda_{j}e_{j}$, that is 
$$ Au_{j} = \lambda_{j}u_{j}\, . $$
As a matter of fact one sees that $\lambda_{j}\in {\Bbb R}$, while  $u_{j}$ is an eigenvector of  $A$ and  $U_{0}$ is a diagonalization operator for $A$, which consists in the matrix whose columns are the eigenvectors $u_{j}$.
Actually this procedure allows us to construct all the eigenvectors of $A$ through the minimization of a unique functional defined on the group ${\bf U}(n)$. Also, since $F(U_{0}) = \sum_{j \geq 1}\lambda_{j}\rho_{j}$, one easily sees that different choices in ordering the numbers $\rho_{j}$ yield different ordering of eigenvalues and eigenvectors of $A$: for instance one may check that if the $\rho_{j}$'s are decreasing, that is if $\rho_{j} > \rho_{j+1}$ for $1\leq j \leq n-1$, then one obtains the eigenvalues of $A$ in a non decreasing order, that is  $\lambda_{1} \leq \lambda_{2}\leq\cdots \leq \lambda_{n}$. While if  $\rho_{j} < \rho_{j+1}$ for $1\leq j \leq n-1$, then one obtains the eigenvalues in a non increasing order, that is  $\lambda_{1} \geq \lambda_{2}\geq\cdots \geq\lambda_{n}$. (Had we began by maximizing $J$, the conclusion would be somehow reversed but analogous: if the $\rho_{j}$'s are decreasing then the $\lambda_{j}$'s would be non increasing). 

\bigskip

In the next section of this paper we will show that, for a certain class of self-adjoint operators $A$, the eigenvectors and eigenvalues of $A$ can be obtained through the minimization of the functional
$$ J_{0}(U) := {\rm tr}(DU^*AU) $$
on an appropriate subset of unitary operators $U$: this is precisely stated and proved in section \ref{sec:Linear}. In section \ref{sec:Prelim} we gather a certain number of preliminary results used in section \ref{sec:SP}, after stating the assumptions on the domain $\Omega$ and on the sequence $(\rho_{m})_{m}$,  we prove Theorem \ref{lem:MainThm}, as well as slightly more general variants of the Schr\"odinger--Poisson systems (see Theorem \ref{lem:SolSP} in  section \ref{sec:SP}).
In section \ref{sec:Remarks} we shall discuss some generalizations and state a few remarks about the results presented here.


\section{Global determination of eigenvectors and eigenvalues}\label{sec:Linear}

In this section we consider an infinite dimensional, separable, complex Hilbert space $H$ whose scalar product is denoted by $(\cdot|\cdot)$ and its norm by $\|\cdot\|$. We shall make the following assumptions:

\begin{hypothesis}\label{lem:HypOnA} 
We assume that $(A,D(A))$ is a densely defined, selfadjoint positive oper\-ator acting on the Hilbert space $H$, and that the domain $D(A)$ equipped with its graph norm is compactly imbedded in $H$, so that $A$ has a compact resolvent and $A$ possesses a sequence of eigenvalues $(\mu_{j})_{j \geq1}$ such that $0 \leq \mu_{j} < \mu_{j+1}$, each eigenvalue having finite multiplicity $m_{j} \geq 1$, and $\mu_{j} \to + \infty$ as $ j \to \infty$, $H$ being infinite dimensional.
\end{hypothesis}

\noindent We denote by $D(A^{1/2})$ the domain of $A^{1/2}$, that is the subspace of $H$ obtained upon the completion of $D(A)$ with the scalar product $(u,v) \mapsto (u|v) + (Au|v)$, and we recall that $D(A^{1/2})$ is dense in $H$. Hence we can introduce the next assumption:

\begin{hypothesis}\label{lem:HypOnej} 
We consider a fixed Hilbert basis of $H$, denoted by $(e_{j})_{j \geq1}$, such that $e_{j} \in D(A^{1/2})$ for each $j \geq 1$.
\end{hypothesis}

\noindent With the Hilbert basis $(e_{j})_{j\geq1}$ given by hypothesis~\lemref{lem:HypOnej}, we consider a sequence $(\rho_{j})_{j \geq 1}$ of real numbers such that
\begin{equation}\label{eq:CondRhoj}
\rho_{j} > 0, \qquad \sum_{j\geq 1}\rho_{j} \, \|e_{j}\|_{D(A^{1/2})}^2 < \infty,
\end{equation}
and we denote by $D$ the diagonal operator defined by 
\begin{equation}\label{eq:DefD}
De_{j} := \rho_{j}e_{j},\qquad\mbox{for }\; j\geq 1.
\end{equation} 
Note that since $H$ is infinite dimensional and $A$ has a compact resolvent, while $e_{j} \in D(A^{1/2})$, we have $\|e_{j}\|_{D(A^{1/2})} \to \infty$ as $j \to \infty$. Indeed, otherwise, the sequence $(e_{j})_{j}$ would be bounded in $D(A^{1/2})$, and the imbedding $D(A^{1/2}) \subset H$ being compact, one would extract a subsequence $(e_{j_{k}})_{k\geq1}$ such that $e_{j_{k}} \weak f$ in $D(A^{1/2})$ and $e_{j_{k}} \to f$ strongly in $H$; in particular $\|f\| = 1$,  since $(e_{j})_{j}$ is a Hilbert basis of $H$. But we have also $e_{j} \weak 0$ in $H$, and thus we should have $f=0$. This contradiction shows that $(e_{j})_{j}$ cannot contain any bounded sequence in $D(A^{1/2})$. As a consequence we have $\rho_{j} \to 0$ and $D$ is a compact operator. 

Next we shall consider unitary operators $U : H \dans H$ which satisfy the following condition (this expresses the fact that the operator $DU^*AU$ is of trace class, see M.~Reed \& B.~Simon 
\cite{ReedSimon:1}, volume 1, section VI.6)
\begin{equation}
\begin{cases}\label{eq:CondU-1}
U^*U=UU^* = I, \quad
Ue_{j}\in D(A^{1/2}) \;\mbox{for }\, j \geq1, \\ 
\sum_{j\geq1} \rho_{j} (U^*AUe_{j}|e_{j}) < \infty.
\end{cases}
\end{equation}
and we define the set ${\Bbb S}$ through
\begin{equation}\label{eq:DefS-1}
{\Bbb S} := \big\{ U : H \dans H \; ; \; U \, \mbox{satisfies \eqref{eq:CondU-1}} \big\}.
\end{equation}

\begin{remark}\label{lem:SNotEmpty} 
Let us point out that such operators $U$ exist, that is ${\Bbb S}$ is not empty: indeed, not only for $t\in{\Bbb R}$ we have ${\rm e}^{{\rm i}tA} \in {\Bbb S}$, but also for any $\lambda>0$, the operator $U_{\lambda}$ (the so-called Cayley transform of $\lambda A$, see K.~Yosida \cite{Yosida:Book}) defined by
$$ U_{\lambda} := (I + {\rm i}\lambda A)(I - {\rm i}\lambda A)^{-1} $$
 is a bounded operator on $H$ and one checks easily that 
$$ U_{\lambda}^* = (I-{\rm i}\lambda A)(I + {\rm i}\lambda A)^{-1},$$
so that $U_{\lambda}^*U_{\lambda} = I$. Moreover, for any $f \in D(A^{1/2})$ we have $(I+{\rm i}\lambda A)^{-1}f \in D(A^{3/2})$, and thus $U_{\lambda}f \in D(A^{1/2})$. As a matter of fact, not only $U_{\lambda}$ is a unitary operator on $H$, but one has also $\|U_{\lambda}f\|_{D(A^{1/2})} = \|f\|_{D(A^{1/2})}$. Therefore, since the sequence $(\rho_{j})_{j}$ satisfies \eqref{eq:CondRhoj}, one sees that $U_{\lambda}$ satisfies \eqref{eq:CondU-1} and $U_{\lambda} \in {\Bbb S}$. 
\qed
\end{remark}
\bigskip
For a unitary operator $U : H \dans H$ satisfying \eqref{eq:CondU-1}, we define $J_{0}(U)$ by
\begin{equation}\label{eq:DefJ0-1}
J_{0}(U) := {\rm tr}(DU^*AU) := \sum_{j\geq 1}\rho_{j}(U^*AUe_{j}|e_{j})
\end{equation}
The following result concerns eigenvectors of $A$:

\begin{theorem}\label{lem:EigenA}
Assume that the hypotheses~\lemref{lem:HypOnA} and \lemref{lem:HypOnej}, as well as condition \eqref{eq:CondRhoj} are satisfied. Then the functional $J_{0}$ defined in \eqref{eq:DefJ0-1} achieves its minimum on ${\Bbb S}$ defined by \eqref{eq:DefS-1}. Then there exists ${\widehat U}_{0} \in {\Bbb S}$ such that
$$ J_{0}({\widehat U}_{0}) = \min_{U \in {\Bbb S}} J_{0}(U), $$
and for each $j \geq 1$, the vector $\phi_{j} := {\widehat U}_{0}e_{j}$ is an eigenvector of $A$ corresponding to the eigenvalue $\lambda_{j} := (A\phi_{j}|\phi_{j})$.
\end{theorem}

\begin{remark}
It is clear that the eigenvalues $(\lambda_{j})_{j\geq1}$ are independent of the choice of the sequence $(\rho_{j})_{j\geq1}$. However, as pointed out in the introduction, if one assumes that the sequence $(\rho_{j})_{j\geq 1}$ is decreasing, then the eigenvalues $(\lambda_{j})_{j\geq1}$ are ordered in a non decreasing order, that is $\lambda_{j} \leq \lambda_{j+1}$ for $j \geq1$. Thus for different choices of the sequence $(\rho_{j})_{j\geq 1}$ one may obtain different diagonalization operators ${\widehat U}_{0}$ for $A$. 
\qed
\end{remark}
We split the proof of this result into a couple of lemmas.

\begin{lemma}\label{lem:MinAchieved} 
The functional $J_{0}$ achieves its minimum on ${\Bbb S}$ at a certain $U_{0}\in {\Bbb S}$.
\end{lemma}

\proof Indeed consider the infimum $\alpha := \inf_{U \in {\Bbb S}} J_{0}(U)$. 
Since ${\Bbb S} \neq \emptyset$, we have $0 \leq \alpha¬†< \infty$. Consider a minimizing sequence $(U_{n})_{n\geq 1} \in {\Bbb S}$, such that for instance $\alpha \leq J(U_{n}) \leq \alpha + 1/n$. Then for each fixed $j \geq 1$, setting $u_{j}^n := U_{n}e_{j}$, we have for all $n \geq1$
$$ \|u_{j}^n\|_{D(A^{1/2})}^2 = 1 + (Au_{j}^n|u_{j}^n) \leq 1 + {\alpha + 1 \over \rho_{j}}. $$
Thus, since the inclusion $D(A^{1/2}) \subset H$ is compact, upon extracting subsequences through Cantor's diagonal scheme, and denoting again this diagonal subsequence by $(u_{j}^n)_{n}$, we may assume that for all $j \geq 1$ there exist a family $(u_{j})_{j}$ such that for $j \geq 1$ fixed
$$ u_{j}^n \weak u_{j}\quad\mbox{weakly in }\, D(A^{1/2}),\qquad u_{j}^n \to u_{j}\quad\mbox{strongly in }\, H $$
as $n \to \infty$. Setting $U_{0}e_{j} := u_{j}$, one checks easily that $U_{0}$ can be extended by linearity to the subspace ${\rm span}\{e_{j} \; ; \; j \geq 1\}$, and that for $f \in {\rm span}\{e_{j} \; ; \; j \geq 1\}$ we have
$$ \|U_{0}f\|^2 = \lim_{n\to \infty}\|U_{n}f\|^2 = \|f\|^2. $$
In other words $U_{0}$ is a unitary operator on (the algebraic) ${\rm span}\{e_{j} \; ; \; j \geq 1\}$, and therefore can be extended as such to the whole space $H$. 
Since for any $m \geq 1$ we have
$$
\sum_{j=1}^m \rho_{j}(Au_{j}|u_{j}) \leq \liminf_{n\to \infty} \sum_{j=1}^m \rho_{j}(Au_{j}^n|u_{j}^n) \leq \liminf_{n\to \infty} J_{0}(U_{n}) = \alpha, $$
upon letting $m \to \infty$ we conclude that $J_{0}(U_{0}) \leq \alpha$. Thus, having $U_{0}^*U_{0} = I$ and  $U_{0}e_{j} \in D(A^{1/2})$ for all $j\geq 1$, and $J_{0}(U_{0}) < \infty$, we have $U_{0}\in {\Bbb S}$ and $J_{0}(U_{0}) = \alpha$. \endproof

\medskip
Next we show that $U_{0}$, given by Lemma~\lemref{lem:MinAchieved} is a diagonalization operator for $A$.

\begin{lemma}\label{lem:ADiag}
Under the assumptions of Theorem \lemref{lem:EigenA}, let $U_{0}$ be given by Lemma~\lemref{lem:MinAchieved},  and set $u_{j}:= U_{0}e_{j}$, for $j \geq 1$.
\renewcommand{\labelenumi}{(\roman{enumi})}
\begin{enumerate}
\item Assume that $k \geq 1$ is such that
\begin{equation}\label{eq:RhojDist}
\rho_{\ell}\neq \rho_{k}\quad\mbox{for }\, \ell \neq k.
\end{equation}
Then there exist $\lambda_{k}\in {\Bbb R}_{+}$ such that $Au_{k} = \lambda_{k}u_{k}$.
\item Assume that $k \geq 1$ is such that for some $m \geq 2$
\begin{equation}
\begin{cases}\label{eq:RhojMult}
\rho_{k} = \rho_{k+\ell}\quad \mbox{for }\, 0 \leq \ell \leq m-1,\\ 
\rho_{k} \neq \rho_{n} \quad \mbox{for }\, n \not\in \left\{k+\ell \; ; \; 0 \leq \ell \leq m-1 \right\}.
\end{cases}
\end{equation}
Then there exists a unitary transformation $U_{k}$ of the $m$-dimensional space $H_{k} := {\rm span}\left\{ U_{0}e_{k+\ell} \; ; \; 0 \leq \ell \leq m-1 \right\}$ such that if 
$${\widehat u}_{k+\ell} := U_{k}U_{0}e_{k+\ell},$$
then there exists $\lambda_{k+\ell} \in {\Bbb R}_{+}$ such that $A{\widehat u}_{k+\ell} = \lambda_{k+\ell}{\widehat u}_{k+\ell}$ for $0 \leq \ell \leq m-1$.
\end{enumerate} 
\end{lemma}

\proof First let $M : H \dans H$ be a bounded skewadjoint operator such that $M : D(A^{1/2}) \dans D(A^{1/2})$ is also bounded. Indeed such operators do exist (consider for instance ${\rm i}(I+\lambda A)^{-1}$ for $\lambda>0$). Setting $U(t) := \exp(-tM)U_{0}$ for $t \in {\Bbb R}$, one checks easily that, since $M^* = -M$, one has $U(t) \in {\Bbb S}$ for all $t$, and thus the function $g(t) := J_{0}(U(t))$ is well defined, is of class $C^1$ and achieves its minimum at $t=0$. However since
$$ g(t) = {\rm tr}(DU_{0}^* \exp(tM)A\exp(-tM)U_{0}), $$
one concludes that 
\begin{equation}\label{eq:M:Skew}
g'(0) = {\rm tr}(DU_{0}^* MAU_{0}) - {\rm tr}(DU_{0}^*AMU_{0}) = 0
\end{equation}
for all bounded operators $M : H \dans H$ such that $M^*= - M$ and $M$ is also bounded from $D(A^{1/2})$ into itself. 
In the same way, if we consider a bounded operator $L: H \dans H$ such that $L = L^*$ and $L$ is also bounded from $D(A^{1/2})$ into itself, upon setting $M := {\rm i}L$, we conclude that \eqref{eq:M:Skew} yields
\begin{equation}\label{eq:L:Auto}
{\rm tr}(DU_{0}^* LAU_{0}) = {\rm tr}(DU_{0}^*ALU_{0}),
\end{equation}
for all such operators $L$.

Note that the above relation \eqref{eq:M:Skew} yields that
\begin{equation*}
\begin{aligned}
& \sum_{j \geq 1} \rho_{j}(U_{0}^*MAU_{0}e_{j}|e_{j})  - 
\sum_{j \geq 1} \rho_{j}(U_{0}^*AMU_{0}e_{j}|e_{j}) = 0\\
& \sum_{j \geq 1} \rho_{j}(MAu_{j}|u_{j}) -
\sum_{j \geq 1} \rho_{j}(AMu_{j}|u_{j}) = 0 ,
\end{aligned}
\end{equation*}
that is, since $M^* = -M$,
\begin{equation}\label{eq:M:Rel}
- \sum_{j \geq 1} \rho_{j}(Au_{j}|Mu_{j}) -
\sum_{j \geq 1} \rho_{j}(AMu_{j}|u_{j}) = 0 \iff
{\rm Re}\sum_{j \geq 1} \rho_{j}(Au_{j}|Mu_{j}) = 0.
\end{equation}
Analogously using \eqref{eq:L:Auto} one obtains in the same way
\begin{equation}\label{eq:L:Rel}
\sum_{j \geq 1} \rho_{j}(Au_{j}|Lu_{j}) = 
\sum_{j \geq 1} \rho_{j}(ALu_{j}|u_{j}) \iff
{\rm Im}\sum_{j \geq 1} \rho_{j}(Au_{j}|Lu_{j}) = 0.
\end{equation}

At this point, in a first step, assume that the integer $k$ is such that condition \eqref{eq:RhojDist} is fulfilled.
Consider an integer $n \neq k$, so that $\rho_{n} \neq \rho_{k}$, and define the operators $M$ and $L$ in the following way
\begin{equation}\label{eq:ChoiceML}
\begin{cases}
Mu_{k} := u_{n}, \quad Mu_{n} := - u_{k}, \\
Lu_{k} := u_{n}, \quad Lu_{n} := u_{k}, \\ 
Lu_{j} = Mu_{j} = 0 \; \mbox{ for }\, j \not\in \{k,n\}.
\end{cases}
\end{equation}

Clearly $M$ and $L$ satisfy the required conditions above, and using \eqref{eq:M:Rel}, with our choice of the operator $M$, we get $(\rho_{n} -\rho_{k}){\rm Re}(Au_{k}|u_{n}) = 0$, that is, since $\rho_{k} - \rho_{n} \neq 0$,
$$ {\rm Re}(Au_{k}|u_{n}) = 0. $$
Upon using \eqref{eq:L:Rel}, with our above choice of the operator $L$ and the fact that $\rho_{n}-\rho_{k} \neq 0$, analogously we have that
$$ {\rm Im}(Au_{k}|u_{n}) = 0. $$
So, from the above two relations, we infer that $(Au_{k}|u_{n}) = 0$ for all $n$ such that $\rho_{n}\neq \rho_{k}$, that is 
$$ Au_{k} \in {\rm span}\{u_{n} \; ; \; n\neq k\}^\perp = {\rm span}\{u_{k}\}, $$
where we use the fact that the family $(u_{j})_{j}$ is a Hilbert basis of $H$, being the image of the Hilbert basis $(e_{j})_{j}$ under the unitary operator $U_{0}$. This means that $Au_{k} = \lambda_{k}u_{k}$ for some $\lambda_{k} \in {\Bbb C}$, but since $A$ is a nonnegative self-adjoint operator, as a matter of fact we have $\lambda_{k} \geq 0$. 

Next assume that the integer $k$ is such that the coefficient $\rho_{k}$ has multiplicity $m \geq 2$, that is condition \eqref{eq:RhojMult} is satisfied.
Arguing as above, we consider the following operators $M$ and $L$: for 
$n \not\in \left\{k+j \; ; \; 0 \leq j \leq m-1 \right\}$ and $0 \leq \ell \leq (m-1)$ fixed, set
\begin{equation}\label{eq:ChoiceML-2}
\begin{cases}
Mu_{k+\ell} := u_{n}, \quad Mu_{n} := - u_{k+\ell}, \\ 
Lu_{k+\ell} := u_{n}, \quad Lu_{n} := u_{k+\ell} \\
Lu_{j} = Mu_{j} = 0 \quad \mbox{for all }\, j \not\in \{n,k+\ell\}.
\end{cases}
\end{equation}
Then, proceeding as above, 
we conclude that $(Au_{k+\ell}|u_{n}) = 0$ for all $n \not\in \left\{k+j \; ; \; 0 \leq j \leq m-1 \right\}$, that is:
$$
Au_{k+\ell} \in \left({\rm span}\left\{ u_{n} \; ; \; n \neq k+j, \;\; 0 \leq j \leq m-1 \right\}\right)^\perp $$
that is 
$$ Au_{k+\ell} \in {\rm span}\left\{u_{k+i} \; ; \; 0 \leq i \leq m-1 \right\}. $$
This means that if we set $H_{k} := {\rm span}\left\{u_{k+i} \; ; \; 0 \leq i \leq m-1 \right\}$, then $A : H_{k} \dans H_{k}$ is a self-adjoint operator on the finite dimensional space $H_{k}$. Therefore there exists a unitary operator $U_{k}$, acting on this space, such that if for $0 \leq \ell \leq m-1$ we set ${\widehat u}_{k+\ell} = U_{k}u_{k+\ell} = U_{k}U_{0}e_{k+\ell}$, we have $A{\widehat u}_{k+\ell}=\lambda_{k+\ell}{\widehat u}_{k+\ell}$ for some $\lambda_{k+\ell} \geq 0$.
\endproof
\medskip

As we may see from the above analysis, when all the $\rho_{j}$'s are distinct, then $U_{0}$, any unitary operator which minimizes $J_{0}$, is a diagonalization operator for $A$. However in the general case, when some of the coefficients $\rho_{k}$ have multiplicity $m_{k} \geq 2$, it is possible that one has to impose a unitary transformation $U_{k}$ in the space 
$$ H_{k}:={\rm span}\{U_{0}e_{k+\ell} \; ; \; 0 \leq \ell \leq m_{k} - 1\} $$
in order to have the operator $A$ diagonalized. 
In other words, one may find a unitary operator $U_{k}$ on $H_{k}$ such that if $A_{k} := A_{|H_{k}}$ is the trace of $A$ on $H_{k}$, the operator  $U_{k}^*A_{k}U_{k}$ is diagonal. Thus since $\rho_{k}=\rho_{k+\ell}$ for $0 \leq \ell \leq m-1$, if we denote by ${\widehat U}$ the unitary operator obtained through the composition of all such operators $U_{k}$ and $U_{0}$, one has $J_{0}({\widehat U}) = J_{0}(U_{0})$. More precisely, we can we state the following corollary, which ends the proof of Theorem \lemref{lem:EigenA}:

\begin{corollary}\label{lem:ADiag2} 
Under the assumptions of Theorem \lemref{lem:EigenA}, let $U_{k}$ be given by Lemma~\lemref{lem:ADiag} when $k \geq1$ is such that \eqref{eq:RhojMult} is satisfied. Then the operator ${\widehat U}_{0}$ defined by ${\widehat U}_{0}e_{k} = U_{0}e_{k}$ when $k$ satisfies \eqref{eq:RhojDist}, and 
$$
{\widehat U}_{0}e_{k+\ell} := U_{k}U_{0}e_{k+\ell}, \quad\mbox{for }\, 0 \leq \ell \leq m-1, \quad\mbox{when \eqref{eq:RhojMult} is satisfied},
$$
belongs to ${\Bbb S}$, while $J_{0}({\widehat U}_{0}) = J_{0}(U_{0})$ and ${\widehat U}_{0}^*A{\widehat U}_{0}$ is diagonal. Setting $\phi_{j} := {\widehat U}_{0}e_{j}$ for $j \geq 1$, then there exists $\lambda_{j} \geq 0$ such that $A\phi_{j} = \lambda_{j}\phi_{j}$.
\end{corollary}

\begin{remark}
Regarding the finite dimenional case, after the completion of this work and its submission for publication, R.V. Kohn, in a private communication, pointed to one of the authors (O.K.) that L. Mirsky \cite{Mirsky:1975}, develops a result of J. von Neumann \cite{vonNeumann:1937} stating that 
$$|{\rm tr}(AB) | \leq \sum_{j=1}^n \sigma_{j}(A)\sigma_{j}(B)$$
for two $n\times n$ matrices $A,B$, where $(\sigma_{j}(M))_{1\leq j \leq n}$ denotes the decreasing singular values of a matrix $M$. From this L. Mirsky concludes that actually one has
$$\sup_{U,V}|{\rm tr}(BUAV)| = \sum_{j=1}^n \sigma_{j}(A)\sigma_{j}(B),$$
which is another result due to J. von Neumann.
In this respect, in the finite dimensional case our result can be compared to the above result in the particular case of self-adjoint matrices, and moreover the approach given here characterizes the diagonalization matrix by a variational method.
\qed 
\end{remark}


\section{Preliminary results for Schr\"odinger--Poisson system}\label{sec:Prelim}

In this section we prove an existence result regarding the system of equations \eqref{eq:1}--\eqref{eq:ortho}. We shall assume that 
\begin{equation}\label{eq:HypDomain}
\Omega \subset {\Bbb R}^N \; \mbox{ is a bounded domain and that }\; N \leq 6,
\end{equation}
and we endow the (complex) space $L^2(\Omega)$ with its scalar product denoted by $(\cdot|\cdot)$ and its norm $\|\cdot\|$.
Let ${\widetilde V}_{0}$ be a real valued potential such that
\begin{equation}\label{eq:HypPotential}
{\widetilde V}_{0}^+ \in L^1(\Omega),\qquad {\widetilde V}_{0}^- \in L^{p_{0}}(\Omega) \; \mbox{ for some }\; p_{0} > {N \over 2}\;
\mbox{ and }\; p_{0} \geq 1.
\end{equation}
Then we define an unbounded operator $(A,D(A))$ by setting 
\begin{equation}\label{eq:DefOpA}
D(A) := \left\{ u \in H^1_{0}(\Omega) \; ; \; -\Delta u + {\widetilde V}_{0}u \in L^2(\Omega)\right\},\qquad Au := -\Delta u + {\widetilde V}_{0}u.
\end{equation}
This operator is self-adjoint, has a compact resolvent, and there exists a Hilbert basis of eigensystem denoted by $(\lambda_{m},\phi_{m})_{m \geq 1}$, that is (here $\delta_{mn}$ being the Kronecker symbol)
$$ -\Delta \phi_{m} + {\widetilde V}_{0}\phi_{m} = \lambda_{m}\phi_{m},\qquad 
\phi_{m} \in H^1_{0}(\Omega), \qquad 
\int_{\Omega}\phi_{m}(x)\phi_{n}(x)\,dx = \delta_{mn}. $$
It is well-known that by Weyl's theorem there exist two positive constants $c_{1}$, $c_{2}$, depending on $\Omega$ and ${\widetilde V}_{0}$, and $m_{0} \geq 1$ large enough such that for all integers $m \geq m_{0}$ one has
$$ c_{1} m^{2/N} \leq \lambda_{m} \leq c_{2} m^{2/N}. $$
(See for instance \cite{OK:Book}, chapter 5, \S~3, where the case of Neumann boundary conditions is also treated). For this reason, as far as the sequence $(\rho_{m})_{m \geq 1}$ is concerned, in order to ensure the finiteness of the functionals we are going to minimize, we assume that the growth condition \eqref{eq:RhojGrowth} is satisfied.

For a given unitary operator $U : L^2(\Omega) \dans L^2(\Omega)$ such that moreover $U : H^1_{0}(\Omega) \dans H^1_{0}(\Omega)$ is also a bounded operator, we shall denote by $V := V[U]$ the potential defined by the Poisson equation (here $|U\phi_{j}|$ denotes the modulus of the function $U\phi_{j}$)
\begin{equation}\label{eq:DefV}
-\Delta V = \sum_{j \geq 1} \rho_{j} |U\phi_{j}|^2 \quad \mbox{in }\; \Omega, \qquad V \in H^1_{0}(\Omega).
\end{equation}

\begin{remark}\label{rem:CondN}
We wish to explain here the limitation $N \leq 6$ in Theorem \ref{lem:MainThm}. Indeed, thanks to the Sobolev imbedding theorem we have $H^1_{0}(\Omega) \subset L^{2^*}(\Omega)$, where $2^* = 2N/(N-2)$ when $N \geq 3$, while $2^*$ can be any finite exponent if $N = 2$, and $2^* = \infty$ if $N = 1$. This means that the right hand side of the above equation \eqref{eq:DefV}, that is
$$f := \sum_{j \geq 1} \rho_{j} |U\phi_{j}|^2 ,$$
belongs to $L^q(\Omega)$ where $q := N/(N-2)$ if $N \geq 3$, or $q < \infty$ arbitrary if $N = 2$, or $q = \infty$ if $N = 1$. Since $q \geq (2^*)' := 2N/(N+2)$ if and only if $N \leq 6$, we conclude that for such $N$'s we have $f \in H^{-1}(\Omega) \cap L^{2N/(N+2)}(\Omega)$: this ensures that equation \eqref{eq:DefV} has a unique solution $V := V[U] \in H^1_{0}(\Omega)$ when $N \leq 6$ and thus the functional
$$U \mapsto J_{1}(U) := {1 \over 2}\int_{\Omega}|\nabla V[U]|^2dx = 
{1 \over 2}\sum_{j \geq 1}\rho_{j}\int_{\Omega}V[U]|U\phi_{j}|^2 dx$$
is well defined (this functional $J_{1}$ is used in the minimization procedure, see below). Moreover while for $N \leq 5$ one may show that $J_{1}$ is, in an appropriate sense, weakly sequentially continuous, the case $N = 6$ is a limit case and we are only able to prove that $J_{1}$ is weakly sequentially lower semi-continuous. Thus for $N \leq 6$ we shall show that $J_{1}$ is weakly sequentially lower semi-continuous, and this allows us to proceed with our minimization procedure.

Finally, when $N \geq 7$ the functional $J_{1}$ is not well defined for an arbitrary unitary operator $U : L^2(\Omega \dans L^2(\Omega)$ which is also a bounded operator on $H^1_{0}(\Omega)$, and the method we use here has to be modified by considering other classes of operators $U$. 
\qed 
\end{remark}

\begin{remark}\label{rem:EstimV}
Observe that if $V[U]$ satisfies \eqref{eq:DefV}, then by the maximum principle we have $V[U] > 0$ in $\Omega$.
Note also that when $N \leq 3$, we have  $f \in L^q(\Omega)$ for some $q > N/2$. Then a classical regularity result (see for instance \cite{OK-Ultra}, or G. Stampacchia \cite{Stampacchia}) states that there exists a constant $c >0$ such that if $V \in H^1_{0}(\Omega)$ satisfies $-\Delta V = f $ and $f \in L^q(\Omega)$, then
\begin{equation}\label{eq:EstimV}
\|V\|_{\infty} \leq c \, \|f\|_{q}\, .
\end{equation}
Therefore, if $V := V[U]$ is given by \eqref{eq:DefV} we have $V \in L^\infty(\Omega)$ when $N \leq 3$. 

When $N = 4$, we have $f \in L^2(\Omega)$, and in this case, since $N/2 = 2$, we have that $V \in L^p(\Omega)$ for all $p < \infty$. Moreover for a constant $c(p)$ depending on $\Omega$ we have:
\begin{equation}\label{eq:EstimV2}
\|V\|_{p} \leq c(p) \, \|f\|_{2}\, .
\end{equation}
 Finally, when $5 \leq N \leq 6$, we have $f \in L^q(\Omega)$ with 
 $1 < q = N/(N-2) < N/2$: using a regularity result, we have $V \in L^p(\Omega)$ where $p$ is given by
$${1 \over p} = {1 \over q} - {2 \over N},$$
that is $p := N/(N-4)$ (again cf.~\cite{OK-Ultra}, or G. Stampacchia \cite{Stampacchia}). In this case there exists a constant $c > 0$ depending on $\Omega$ such that
\begin{equation}\label{eq:EstimV3}
\|V\|_{N/(N-4)} \leq c \, \|f\|_{N/(N-2)}\, .
\end{equation}
These observations will allow us to show that the functional $J_{1}$ is, in an aproriate sense, weakly sequentially lower semi-continuous (see below Lemma~\lemref{lem:J1Cont}).
\qed
\end{remark}

It is convenient to consider the Sobolev space ${\Bbb H}_{1}$ endowed with the norm $\|\cdot\|_{{\Bbb H}_{1}}$:
$$
{\Bbb H}_{1} := \left\{u \in H^1_{0}(\Omega) \; ; \; \|u\|_{{\Bbb H}_{1}}^2 := \|\nabla u\|^2 + \int_{\Omega}{\widetilde V}_{0}^+(x)|u(x)|^2dx < \infty\right\}.
$$
The imbedding ${\Bbb H}_{1} \subset L^2(\Omega)$ is compact. Note that since the eigenfunctions $\phi_{m}$ belong to $L^\infty(\Omega)$, we have $\phi_{m} \in {\Bbb H}_{1}$.

Regarding the manifold ${\Bbb S}$ defined in \eqref{eq:DefS}, we have to modify it slightly, as we did in section \S~2. More precisely we shall consider unitary operators $U : L^2(\Omega) \dans L^2(\Omega)$ such that
\begin{equation}\label{eq:CondU-2}
U^*U=UU^* = I, \quad
U\phi_{j}\in {\Bbb H}_{1} \;\mbox{ for }\, j \geq1, \quad
\sum_{j\geq1} \rho_{j} (U^*AU\phi_{j}|\phi_{j}) < \infty,
\end{equation}
and we consider the manifold defined by
\begin{equation}\label{eq:DefS-2}
{\Bbb S} := \big\{ U : L^2(\Omega) \dans L^2(\Omega) \; ; \; U \, \mbox{satisfies \eqref{eq:CondU-2}} \big\}.
\end{equation}

We denote by $D$ the diagonal operator acting on $L^2(\Omega)$ defined by $D\phi_{j} = \rho_{j}\phi_{j}$, and for $U \in {\Bbb S}$ we define the functionals $J_{0}$ and $J_{1}$ as follows: 
\begin{equation}\label{eq:DefJ0-bis}
 J_{0}(U) := {\rm tr}(DU^*AU) = \sum_{j \geq 1}\rho_{j} \int_{\Omega}\left(|\nabla U\phi_{j}|^2(x) + {\widetilde V}_{0}(x) |U\phi_{j}|^2(x) \right)dx 
\end{equation}
and 
\begin{equation} \label{eq:DefJ1}
\begin{aligned}
J_{1}(U) &:= {1 \over 2}\int_{\Omega}|\nabla V[U]|^2(x)dx \\
& = \langle -{1 \over 2}\Delta V[U], V[U] \rangle \\
& = {1 \over 2}\sum_{j \geq 1}\rho_{j}(V[U]U\phi_{j}|U\phi_{j}) = {1 \over 2} {\rm tr}(DU^*V[U]U), \\
\end{aligned}
\end{equation}
where, in the last equality of \eqref{eq:DefJ1}, by an abuse of notation, we denote by $V[U]$ the (linear) multiplication operator $f\mapsto V[U]f$. Since we assume $N \leq 6$, as explained in the above Remark \ref{rem:CondN}, the functional $J_{1}$ is well defined on ${\Bbb S}$.

Note that here the potential ${\widetilde V}_{0}$ may have a negative part, so at some point we will need to ensure that the functional $J_{0}$ is bounded below,  that it is {\it coercive\/} in some sense. More precisely we have:

\begin{lemma}\label{lem:J0Below}
There exists $C \geq 0$ such that for any $U \in {\Bbb S}$ one has 
$$
J_{0}(U) \geq {1 \over 2}\sum_{j \geq 1}\rho_{j}\int_{\Omega}\left(|\nabla U\phi_{j}|^2 + 2 {\widetilde V}_{0}^+|U\phi_{j}|^2\right)dx  - C.
$$
\end{lemma}

\proof Assume that $N \geq 3$ (the case $N \leq 2$ can be handled in a similar way). For $t > 0$ and $u\in H^1_{0}(\Omega)$ such that $\|u\| = 1$ we have
\begin{equation*}
\begin{aligned}
\int_{\Omega}{\widetilde V}_{0}^- |u|^2 dx & = 
	\int_{[{\widetilde V}_{0}^- > t]}{\widetilde V}_{0}^-|u|^2 dx
	+ \int_{[{\widetilde V}_{0}^- \leq t]}{\widetilde V}_{0}^-|u|^2 dx \\
& \leq \int_{\Omega}1_{[{\widetilde V}_{0}^- > t]}{\widetilde V}_{0}^-|u|^2 dx
	+ t  \int_{\Omega}|u|^2 dx \\
& \leq \|1_{[{\widetilde V}_{0}^- > t]}{\widetilde V}_{0}^-\|_{L^{N/2}}\, \|u\|_{L^{2^*}}^2 + t \\
& \leq C_{1}(N)\, {\rm meas}([{\widetilde V}_{0}^- > t])^\theta \,\|{\widetilde V}_{0}^-\|_{L^{p_{0}}} \, \|\nabla u\|^2 + t \\
\end{aligned}
\end{equation*}
where we have used H\"older's inequality twice (once with $N/(N-2)$ and $(N/(N-2))' = N/2$, once with  $p_{0}$ and $N/2$, where $\theta = (2/N) - (1/p_{0}) = (2/3) - (1/p_{0})> 0$). We used also Sobolev's inequality $\|u\|_{2^*} \leq C \|\nabla u\|$. Now, since ${\widetilde V}_{0}^- \in L^{p_{0}}(\Omega)$, we know that ${\rm meas}([{\widetilde V}_{0}^- > t]) \to 0$ as $t \to +\infty$. We choose $t >0$ large enough to ensure that 
$$
C_{1}(N)\, {\rm meas}([{\widetilde V}_{0}^- > t])^\theta\, \|{\widetilde V}_{0}^-\|_{L^{p_{0}}} \leq {1 \over 2}.
$$
Then we have for all $u \in {\Bbb H}_{1}$
$$
\int_{\Omega}|\nabla u|^2dx + \int_{\Omega} {\widetilde V}_{0}|u|^2dx \geq {1 \over 2 } \|\nabla u\|^2 + \int_{\Omega} {\widetilde V}_{0}^+|u|^2dx - t.
$$
Applying this to $u:=U\phi_{j}$, multiplying by $\rho_{j} > 0$ and calculating the sum over $j$ yields the inequality claimed by our lemma, with $C := t\sum_{j\geq 1}\rho_{j}$.
\endproof
\medskip

It is well known that the fact that the functional $u \mapsto \|\nabla u\|^2$ is weakly sequentially lower semi-continuous (l.s.c.) on $H^1_{0}(\Omega)$ plays a crucial role in many minimization problems. Regarding the functional $J_{0}$ we need an analogous property which is stated below:

\begin{lemma}\label{lem:J0sci}
The functional $J_{0}$ is ``weakly sequentially lower semi-continuous'' in the following sense:
let $(U_{n})_{n \geq 1}$ be a sequence in ${\Bbb S}$ such that for some $R>0$ and all $n \geq 1$ one has $J_{0}(U_{n}) \leq R$. Then there exists a subsequence $(U_{n_{k}})_{k}$ such that for any fixed $j \geq 1$ one has $U_{n_{k}}\phi_{j}\weak u_{j}$ in ${\Bbb H}_{1}$ as $k \to +\infty$, and if we set $u^{n_{k}}_{j} := U_{n_{k}}\phi_{j}$ and we define a linear operator $U$ by setting $U\phi_{j} := u_{j}$ we have $u^{n_{k}}_{j} \to u_{j}$ strongly in $L^2(\Omega)$ and almost everywhere on $\Omega$, and
$$
\lim_{k\to\infty}\int_{\Omega}{\widetilde V}_{0}^-|u^{n_{k}}_{j}|^2dx = \int_{\Omega}{\widetilde V}_{0}^-|u_{j}|^2dx, \quad U \in {\Bbb S}, \quad J_{0}(U) \leq \liminf_{k \to \infty}J_{0}(U_{n_{k}}).
$$
\end{lemma}

\proof Assume $N \geq 3$ (the case $N \leq 2$ being analgous). Thanks to Lemma~\lemref{lem:J0Below}, we know that 
$$
\sum_{j \geq 1}\rho_{j}\int_{\Omega}\left(|\nabla U_{n}\phi_{j}|^2 + 2 {\widetilde V}_{0}^+|U_{n}\phi_{j}|^2\right)dx  \leq 2R + 2C =: C_{1}.
$$
This implies that for each $j \geq 1$ fixed the sequence $(u^n_{j})_{n} := (U_{n}\phi_{j})_{n}$ is bounded in ${\Bbb H}_{1}$, more precisely
$\|u^n_{j}\|_{{\Bbb H}_{1}}^2 \leq C/\rho_{j}$. By using Cantor's diagonal scheme and the compactness of the imbedding ${\Bbb H}_{1} \subset L^2(\Omega)$, we may extract a subsequence denoted by $(u^{n_{k}}_{j})_{k\geq 1}$ such that
$$
\begin{cases}
u^{n_{k}}_{j} \weak u_{j}\quad\mbox{weakly in }\, {\Bbb H}_{1},\\
 u^{n_{k}}_{j} \to u_{j}\quad\mbox{strongly in }\, L^2(\Omega),\\
u^{n_{k}}_{j} \to u_{j}\quad\mbox{a.e. in }\, \Omega,
\end{cases}
$$
as $k \to \infty$.
For any $m \geq1$ fixed, we have
$$
\sum_{j=1}^{m}\rho_{j}\|u_{j}\|_{{\Bbb H}_{1}}^2 \leq
\liminf_{k \to \infty}\sum_{j=1}^{m}\rho_{j}\|u^{n_{k}}_{j}\|_{{\Bbb H}_{1}}^2 \leq 
\liminf_{k \to \infty}\sum_{j=1}^{\infty}\rho_{j}\|u^{n_{k}}_{j}\|_{{\Bbb H}_{1}}^2 \leq C\, ,
$$
and finally
\begin{equation}\label{eq:EstimujH1}
\sum_{j=1}^{\infty}\rho_{j}\|u_{j}\|_{{\Bbb H}_{1}}^2 \leq
\liminf_{k \to \infty}\sum_{j=1}^{\infty}\rho_{j}\|u^{n_{k}}_{j}\|_{{\Bbb H}_{1}}^2 \leq C.
\end{equation}
Setting $U\phi_{j} := u_{j}$, one checks easily that $U$ can be extended by linearity to the subspace ${\rm span}\{\phi_{j} \; ; \; j \geq 1\}$, and that for $f \in {\rm span}\{\phi_{j} \; ; \; j \geq 1\}$ we have
$$
\|Uf\|^2 = \lim_{n\to \infty}\|U_{n}f\|^2 = \|f\|^2.
$$
In other words $U$ is a unitary operator on (the algebraic) ${\rm span}\{\phi_{j} \; ; \; j \geq 1\}$, and therefore can be extended as such to the whole space $L^2(\Omega)$. Then \eqref{eq:EstimujH1} shows that $U\in {\Bbb S}$.

We note also that in particular we have 
\begin{equation}\label{eq:EstimV0uj}
\sum_{j\geq 1}\rho_{j}\int_{\Omega}{\widetilde V}_{0}^-(x)|u_{j}(x)|^2dx < \infty.
\end{equation}
 Since we assume $N \geq 3$, the strong convergence of $u^{n_{k}}_{j} \to u_{j}$ in $L^2(\Omega)$ implies (through H\"older's inequality, or interpolation between $L^2(\Omega)$ and $L^{2N/(N-2)}(\Omega)$) that for any fixed $p < N/(N-2)$, and any $j \geq 1$ we have that $u^{n_{k}}_{j} \to u_{j}$ strongly in $L^{2p}(\Omega)$ and a.e.{} in $\Omega$, and thus $|u^{n_{k}}_{j}|^2 \to |u_{j}|^2$ strongly in $L^p(\Omega)$. Since ${\widetilde V}_{0}^- \in L^{p_{0}}(\Omega)$ and $p_{0} > N/2$, taking $p := p_{0}' = p_{0}/(p_{0}-1)$, so that $p < N/(N-2)$, we conclude first that
$$
\lim_{k\to\infty}\int_{\Omega}{\widetilde V}_{0}^-|u^{n_{k}}_{j}|^2dx = \int_{\Omega}{\widetilde V}_{0}^-|u_{j}|^2dx,
$$
and then thanks to \eqref{eq:EstimV0uj} and the monotone convergence theorem,
$$
\lim_{k\to\infty}\sum_{j \geq 1}\rho_{j}\int_{\Omega}{\widetilde V}_{0}^-|u^{n_{k}}_{j}|^2dx = \sum_{j \geq 1}\rho_{j}\int_{\Omega}{\widetilde V}_{0}^-|u_{j}|^2dx.
$$
From this and \eqref{eq:EstimujH1} it is clear that
$$
\begin{aligned}
\sum_{j=1}^{\infty}\rho_{j}\int_{\Omega}\left(|\nabla u_{j}|^2 +
{\widetilde V}_{0}^+(x)|u_{j}|\right)dx & \leq 
\liminf_{k\to \infty} J_{0}(U_{n_{k}}) \\
& \qquad \quad+
\lim_{k\to\infty}\sum_{j \geq 1}\rho_{j}\int_{\Omega}{\widetilde V}_{0}^-|u^{n_{k}}_{j}|^2dx \\
& \leq 
\liminf_{k\to \infty}J_{0}(U_{n_{k}}) \\
& \qquad\quad +
\sum_{j \geq 1}\rho_{j}\int_{\Omega}{\widetilde V}_{0}^-|u_{j}|^2dx, \\
\end{aligned}
$$
which means that $J_{0}(U) \leq \liminf_{k\to\infty}J_{0}(U_{n_{k}})$, as claimed.
\endproof
\medskip

Regarding the functional $J_{1}$ we have the following result:
\begin{lemma}\label{lem:J1Cont} 
Let $1 \leq N \leq 6$.
The functional $J_{1}$ is ``weakly sequentially continuous'' in the following sense:
let $(U_{n})_{n \geq 1}$ be a sequence in ${\Bbb S}$ such that for some $R>0$ and all $n \geq 1$ one has
$$ \sum_{j \geq 1}\rho_{j}\|U_{n}\phi_{j}\|^2_{{\Bbb H}_{1}} \leq R, $$
and such that for any fixed $j \geq 1$ one has 
$$U_{n}\phi_{j}\weak u_{j}\,\mbox{ in }\, {\Bbb H}_{1}, \quad
U_{n}\phi_{j} \to u_{j} \,\mbox{ in }\, L^2(\Omega) \, \mbox{ and a.e.\ on }\,\Omega$$
 as $n \to +\infty$. If we set $u^{n}_{j} := U_{n}\phi_{j}$ and we define a linear operator $U$ by setting $U\phi_{j} := u_{j}$ we have
$$J_{1}(U) \leq \liminf_{n\to\infty} J_{1}(U_{n}).$$
\end{lemma}

\proof This property is due to the positivity of the Green function associated to the Dirichlet problem. Indeed it is well-known that for $f\in L^p(\Omega)$ where $p := 2N/(N+2)$ if $N \geq 3$, or $p > 1$ if $N=2$, or $p=1$ if $N=1$, the solution of the Dirichlet problem
$$-\Delta w = f\quad\mbox{in }\,\Omega,\qquad w \in H^1_{0}(\Omega)$$
may be represented with the Green kernel
$$w(x) = \int_{\Omega}K(x,y)\,f(y)\,dy$$
where $K(x,y) > 0$ on $\Omega\times \Omega$. In particular we have $\|\nabla w\|^2 = \int_{\Omega}f(x)\,w(x)\,dx$, and thus
$$\|\nabla w\|^2 = \int_{\Omega\times\Omega} K(x,y)\,f(y)\,f(x)\,dxdy.$$
Applying this to $f := f^n := \sum_{j\geq 1}\rho_{j}|u^n_{j}|^2$ and $w := w^n := V[U_{n}]$, we see that 
$$J_{1}(U_{n}) = {1 \over 2}\sum_{j,k\geq 1} \rho_{j}\rho_{k} \int_{\Omega\times\Omega} K(x,y)\, |u^n_{j}(x)|^2 |u^n_{k}(y)|^2\,dxdy.$$
Since and $K(x,y) > 0$ and $u^n_{j}(x) \to u_{j}(x)$ a.e.\ on $\Omega$, for any fixed $m \geq 1$, by Fatou's lemma we have
\begin{align*}
&{1 \over 2}\sum_{j,k = 1}^m \rho_{j}\rho_{k} \int_{\Omega\times\Omega} K(x,y)\, |u_{j}(x)|^2 |u_{k}(y)|^2\,dxdy  \leq \cr
&\qquad\qquad {1 \over 2}\liminf_{n\to\infty} \sum_{j,k = 1}^m \rho_{j}\rho_{k} \int_{\Omega\times\Omega} K(x,y)\, |u^n_{j}(x)|^2 |u^n_{k}(y)|^2\,dxdy \leq \cr
&\qquad {1 \over 2}\liminf_{n\to\infty}\sum_{j,k\geq 1} \rho_{j}\rho_{k} \int_{\Omega\times\Omega} K(x,y)\, |u^n_{j}(x)|^2 |u^n_{k}(y)|^2\,dxdy = \liminf_{n\to\infty} J_{1}(U_{n})
\end{align*}
Finally letting $m\to\infty$, we get $J_{1}(U) \leq \liminf_{n\to\infty}J_{1}(U_{n})$. \endproof

\begin{remark} As a matter of fact when $N \leq 5$, we have $J_{1}(U) = \lim_{n\to\infty}J_{1}(U_{n})$, that is $J_{1}$ is weakly sequentially continuous. Indeed let us assume $3 \leq N \leq 5$, as the case $N \leq 2$ can be easily handled analogously.
Since for a fxed $j \geq 1$ we have $u^{n}_{j} \weak u_{j}$ in ${\Bbb H}_{1}$, and $u^{n}_{j} \to u_{j}$ strongly in $L^2(\Omega)$, and a.e.\ on $\Omega$, as $n \to\infty$, we infer that for $2 < p < 2^* = 2N/(N-2)$ we have $u^n_{j} \to u_{j}$ strongly in $L^p(\Omega)$ and that $|u^n_{j}|^2 \to |u_{j}|^2$ in $L^{p/2}(\Omega)$ as $n \to\infty$.
On the other hand, by H\"older's inequality and Sobolev imbedding theorem, if $v \in H^1_{0}(\Omega)$ and $\|v\|_{1} =1$, for $1 < p/2 < N/(N-2)$ we have
$$\|v^2\|_{p/2} \leq \|v^2\|_{1}^{1-\theta}\, \|v^2\|_{N/(N-2)}^\theta = \|v\|_{2^{*}}^{2\theta} \leq 
c\, \|\nabla v\|^{2\theta},$$
with $0 < \theta < 1$ given by
$${2 \over p} = {1-\theta \over 1} + {\theta(N-2) \over N}. $$
Applying this to $v := |u^n_{j}|^2$ yields
$$\left\||u^n_{j}|^2\right\|_{p/2} \leq c\, \|\nabla u^n_{j}\|^{2\theta}.$$
Now we infer first that for any $m \geq 1$ fixed, upon setting $q_{0} := 1 /\theta >1$ and $q'_{0} = q_{0}/(q_{0}-1) < \infty$ we have
$$
\begin{aligned}
\left\| \sum_{j > m} \rho_{j} |u^n_{j}|^2 \right\|_{p/2} & \leq \sum_{j > m} \rho_{j} \|u^n_{j}\|^2_{p} 
 \leq c\, \sum_{j > m}\rho_{j} \|\nabla u^n_{j}\|^{2\theta} \\
&\hskip-1cm \leq c\, \left(\sum_{j > m}\rho_{j}\right)^{1/q'_{0}} 
\left(\sum_{j > m}\rho_{j} \|\nabla u^n_{j}\|^2 \right)^{1/q_{0}} 
\hskip-.5cm  \leq c\, R^{1/q}\, \left(\sum_{j > m}\rho_{j}\right)^{1/q'}.\\
\end{aligned}
$$
From this it follows that
\begin{equation}\label{eq:Def-fn}
f_{n} := \sum_{j \geq 1} \rho_{j} |u^n_{j}|^2 \to f := \sum_{j \geq 1} \rho_{j} |u_{j}|^2 \quad \mbox{strongly in }\, L^{p/2}(\Omega).
\end{equation}
Since here we are assuming that $N \leq 5$, we have $N/(N-2) > 2N/(N+2) = (2^*)'$ and thus we may fix $p > 2$ such that $2N/(N+2) \leq p/2 < N/(N-2)$: therefore $f_{n} \to f$ in $H^{-1}(\Omega)$, since $L^{p/2}(\Omega) \subset H^{-1}(\Omega)$ with continuous imbedding (the functions $f_{n}$ and $f$ are defined in \eqref{eq:Def-fn}).
Therefore, since $-\Delta V[U_{n}] = f_{n}$ and $V[U_{n}] \in H^1_{0}(\Omega)$, clearly we can deduce that $V[U_{n}] \to V[U]$ strongly in $H^1_{0}(\Omega)$ and finally
$$
\int_{\Omega}|\nabla V[U_{n}]|^2\,dx  \to 
\int_{\Omega}|\nabla V[U]|^2\,dx , 
$$
as claimed.\endproof
\end{remark}

\section{Existence of solutions for the Schr\"odinger--Poisson system}\label{sec:SP}

In this section we solve the following Schr\"odinger--Poisson problem:
\begin{numcases}{}
 - \Delta  u_m +{\widetilde V}_{0}u_{m} + V u_m   = \lambda_{m} u_m  \quad\mbox{in }\, \Omega  \label{eq:SP1} \\
- \Delta V  = \sum_{m=1}^{\infty}\rho_{m}| u_m|^2  \quad\mbox{in }\, \Omega  \label{eq:SP2}
\end{numcases}
and moreover
\begin{numcases}{}
 u_m \in H^1_{0}(\Omega),\; \mbox{ for all } m \geq 1, \quad V \in H^1_{0}(\Omega), \label{eq:SP3} \\
(u_m)_{m\geq1} \; \mbox{ is a Hilbert basis for }\, L^2(\Omega). \label{eq:SP4}
\end{numcases}
The main result of this section is:
\begin{theorem}\label{lem:SolSP}
Assume that the hypotheses~\eqref{eq:HypDomain}--\eqref{eq:HypPotential}, as well as condition \eqref{eq:RhojGrowth} are satisfied. The functionals $J_{0}$ and $J_{1}$ being defined in \eqref{eq:DefJ0-bis}--\eqref{eq:DefJ1}, we set $J(U) := J_{0}(U) + J_{1}(U)$ for $U \in {\Bbb S}$ given by \eqref{eq:DefS-2}. Then $J$ achieves its minimum on ${\Bbb S}$ and there exists ${\widehat U}_{0} \in {\Bbb S}$ such that
$J_{0}({\widehat U}_{0}) = \min_{U \in {\Bbb S}} J_{0}(U)$,
and the family $u_{j} := {\widehat U}_{0}\phi_{j}$ is solution to \eqref{eq:SP1}--\eqref{eq:SP4}. Moreover if $V\in H^1_{0}(\Omega)$ satisfies $-\Delta V = \sum_{m \geq 1}\rho_{m}|u_{m}|^2$, then the eigenvalues  $\lambda_{j}$ satisfy 
$$
\lambda_{j} := (-\Delta u_{j} + {\widetilde V}_{0}u_{j} + V u_{j}|u_{j}) = \int_{\Omega}|\nabla u_{j}|^2dx + \int_{\Omega}\left({\widetilde V}_{0} + V\right)u_{j}(x)^2dx.
$$
\end{theorem}

We split the proof of this theorem into several lemmas. First we show that $J$ achieves indeed its minimum.
\begin{lemma}\label{lem:MinAchievedSP}
The functional $J$ achieves its minimum on ${\Bbb S}$ at a certain $U_{0}\in {\Bbb S}$.
\end{lemma}

\proof Since ${\Bbb S} \neq \emptyset$ and $J_{0}$ is bounded below (see Lemma~\lemref{lem:J0Below}), so is $J$ and the infimum
$$ \alpha := \inf_{U \in {\Bbb S}} J(U) $$
is finite. Consider a minimizing sequence $(U_{n})_{n\geq 1} \in {\Bbb S}$, such that for instance $\alpha \leq J(U_{n}) \leq \alpha + 1/n$. 
In particular $J_{0}(U_{n}) \leq 1 + \alpha$, and thanks to Lemma~\lemref{lem:J0sci}, there exists a subsequence (which denote again by $U_{n}$) such that if we set $u_{j}^n := U_{n}\phi_{j}$ for each fixed $j \geq 1$, we have, for all $2\leq p < 2^*$,
$$
\begin{cases}{}
u^{n}_{j} \weak u_{j}\quad\mbox{weakly in }\, {\Bbb H}_{1},\\ 
u^{n}_{j} \to u_{j}\quad\mbox{strongly in }\, L^p(\Omega),\\
u^{n}_{j} \to u_{j}\quad\mbox{a.e. in }\, \Omega ,
\end{cases}
$$
and the operator $U_{0}$ being defined by $U_{0}\phi_{j}:=u_{j}$, we have $U_{0} \in {\Bbb S}$ and
$$
J_{0}(U_{0}) \leq \liminf_{n \to \infty}J_{0}(U_{n}).
$$
On the other hand, thanks to Lemma~\lemref{lem:J1Cont}, we know that $J_{1}$ is weakly sequentially lower semi-continuous, that is $J_{1}(U_{0}) \leq \liminf_{n\to\infty} J_{1}(U_{n})$. Thus
\begin{equation*}
J(U_{0}) = J_{0}(U_{0}) + J_{1}(U_{0}) \leq \liminf_{n\to \infty} J_{0}(U_{n}) + \liminf_{n\to \infty}J_{1}(U_{n}) \leq \liminf_{n\to\infty} J(U_{n}),
\end{equation*}
that is $J(U_{0}) \leq \alpha$. Since $U_{0}\in {\Bbb S}$, we conclude that $J(U_{0})=\alpha$ that is $J$ achieves its minimum at $U_{0}$. 
\endproof

\begin{lemma}\label{lem:DeriveJ0}
Let $U_{0}$ be given by Lemma~\lemref{lem:MinAchievedSP}, and
let $M : H \dans H$ be a bounded skewadjoint operator such that $M : {\Bbb H}_{1} \dans {\Bbb H}_{1}$ is also bounded. Set $U(t) := \exp(-tM)U_{0}$ for $t \in {\Bbb R}$, and $g_{0}(t) := J_{0}(U(t))$. Then $g_{0}$ is of class $C^1$ and 
\begin{equation} 
g'_{0}(0) = - 2{\rm Re} \sum_{j \geq 1} \rho_{j}(AU_{0}\phi_{j}|MU_{0}\phi_{j}).
\end{equation}
\end{lemma}

\proof First one checks easily that, since $M^* = -M$, one has $U(t) \in {\Bbb S}$ for all $t$, and thus the function 
$$ g_{0}(t) := J_{0}(U(t)) = \sum_{j \geq 1}\rho_{j} (AU(t)\phi_{j}|U(t)\phi_{j}) $$
is well defined and is of class $C^1$. Since $U'(t) := dU(t)/dt = -M\exp(-tM)U_{0}$, one sees that
$$
g'_{0}(t) = \sum_{j \geq 1}\rho_{j} (AU'(t)\phi_{j}|U(t)\phi_{j}) +
\sum_{j \geq 1}\rho_{j} (AU(t)\phi_{j}|U'(t)\phi_{j}),
$$
and finally,
$$
g'_{0}(0) = -\sum_{j \geq 1}\rho_{j} (AMU_{0}\phi_{j}|U_{0}\phi_{j}) -
\sum_{j \geq 1}\rho_{j} (AU_{0}\phi_{j}|MU_{0}\phi_{j}),
$$
which yields our claim since $A^* = A$.
\endproof
\medskip

We have an analogous result concerning the functional $J_{1}$: before showing this, we need to show that the mapping $U \mapsto V[U]$ is smooth.

\begin{lemma}\label{lem:DeriveV}
Let $U_{0}$ be given by Lemma~\lemref{lem:MinAchievedSP}, and
let $M : H \dans H$ be a bounded skewadjoint operator such that $M : {\Bbb H}_{1} \dans {\Bbb H}_{1}$ is also bounded. Set $U(t) := \exp(-tM)U_{0}$ for $t \in {\Bbb R}$. Denoting by $V(t)$ the mapping $t \mapsto V[U(t)]$, then $t \mapsto V(t)$ is of class $C^1$ from ${\Bbb R}$ into $L^\infty(\Omega)\cap H^1_{0}(\Omega)$ and denoting by $W$ the solution of
$$
-\Delta W = - 2{\rm Re} \sum_{j \geq 1} \rho_{j} 
\left(MU_{0}\phi_{j}\right) {\overline {U_{0}\phi_{j}}}, \qquad W\in H^1_{0}(\Omega), 
$$ 
we have $V'(0) = W$.
\end{lemma}

\proof The fact that for all $T>0$ and $t \in [-T,T]$ we have
$$
\|M\exp(-tM)U_{0}\phi_{j}\|_{{\Bbb H}_{1}} \leq \|M\|_{{\cal L}({\Bbb H}_{1})}\exp\left(T\|M\|_{{\cal L}({\Bbb H}_{1})}\right)\, \|U_{0}\phi_{j}\|_{{\Bbb H}_{1}} ,
$$ 
shows that for any $p$ such that  $1 < p/2 < N/(N-2) $ (see the proof of Lemma~\lemref{lem:J1Cont}; here assume that $3 \leq N \leq 6$, the case $N \leq 2$ being treated analogously) the mapping 
$$
t \mapsto \sum_{j \geq 1} \rho_{j} |U(t)\phi_{j}|^2
$$
is of class $C^1$ from $(-T,T) \dans L^{p/2}(\Omega) \subset H^{-1}(\Omega)$, and thus using \eqref{eq:EstimV} or \eqref{eq:EstimV2} or \eqref{eq:EstimV3}, the mapping $t \mapsto V(t):=V[U(t)]$ is of class $C^1$ from $(-T,T)$ into $L^{p_{1}}(\Omega)\cap H^1_{0}(\Omega)$, where $p_{1} := \infty$ if $N \leq 3$, and $p_{1} < \infty$ arbitrary if $N = 4$, and $p_{1} := N/(N-4)$ if $5 \leq N \leq 6$ (cf. Remark \ref{rem:EstimV} above). The calculation of $V'(0)$ is straightforward.

\endproof

Now we can state the following result, which will allow us to characterize $U_{0}$ given by Lemma~\lemref{lem:MinAchievedSP}.

\begin{lemma}\label{lem:DeriveJ1}
Let $U_{0}$ be given by Lemma~\lemref{lem:MinAchievedSP}, and
let $M : H \dans H$ be a bounded skewadjoint operator such that $M : {\Bbb H}_{1} \dans {\Bbb H}_{1}$ is also bounded. Set $U(t) := \exp(-tM)U_{0}$ for $t \in {\Bbb R}$, and $g_{1}(t) := J_{1}(U(t))$. Then $g_{1}$ is of class $C^1$ and 
\begin{equation}
g'_{1}(0) = -2\, {\rm Re}\sum_{j \geq 1}\rho_{j}(V[U_{0}]U_{0}\phi_{j}|MU_{0}\phi_{j}).
\end{equation}
\end{lemma}

\proof Thanks to Lemma~\lemref{lem:DeriveV} one checks easily that the function 
$$
g_{1}(t) := J_{1}(U(t)) = {1 \over 2} \int_{\Omega}
\big|\nabla V[U(t)]\big|^2 dx = {1 \over 2} (-\Delta V[U(t)] | V[U(t)])
$$
is $C^1$, and that denoting by $W:=V'(0)$, we have (with $V_{0}:=V[U_{0}])$
$$
\begin{aligned}
g'_{1}(0) = (-\Delta W|V_{0}) & = -2\, {\rm Re}\sum_{j\geq 1}\rho_{j} 
\int_{\Omega}(MU_{0}\phi_{j}){\overline {U_{0}\phi_{j}}} V_{0}dx \\
& = -2\, {\rm Re}\sum_{j \geq 1}\rho_{j}(MU_{0}\phi_{j}|V_{0}U_{0}\phi_{j}), \\
\end{aligned}
$$
where we use the fact that $V_{0}:=V[U_{0}]$ is real valued. 
\endproof

\medskip
The following result is analogous to Lemma~\lemref{lem:ADiag}: the only difference is that due to the presence of the nonlinear term we have to check that when some $\rho_{k}$ has multiplicity $m \geq 2$, we can still proceed as before.

\begin{lemma}\label{lem:SPDiag}
Under the assumptions of Theorem \lemref{lem:SolSP}, let $U_{0}$ be given by lemma~\lemref{lem:MinAchievedSP},  and set $u_{j}:= U_{0}\phi_{j}$, for $j \geq 1$. Then the conclusions of lemma~\lemref{lem:ADiag} hold.
\end{lemma}

\proof The proof is very much the same as in Lemma~\lemref{lem:ADiag}, so we give only the outline and the changes to be made. With the notations of Lemmas \lemref{lem:DeriveJ0} and \lemref{lem:DeriveJ1}, we set $g(t) := J(U(t)) = g_{0}(t) + g_{1}(t) $. Since $g(0) \leq g(t)$ for all $t\in {\Bbb R}$, we have $g'(0)=0$, that is 
\begin{equation}\label{eq:M:Skew-2}
{\rm Re} \sum_{j \geq 1} \rho_{j}(AU_{0}\phi_{j}|MU_{0}\phi_{j})
+ {\rm Re}\sum_{j \geq 1}\rho_{j}(V[U_{0}]U_{0}\phi_{j}|MU_{0}\phi_{j})
= 0
\end{equation}
for all bounded skew-adjoint operators $M$ such that $M : {\Bbb H}_{1} \dans {\Bbb H}_{1}$ is also bounded. 

In the same way, if we consider a bounded adjoint operator $L$ such that $L: {\Bbb H}_{1} \dans {\Bbb H}_{1}$ is also bounded, we may set $M:= {\rm i}\, L$ and conclude that \eqref{eq:M:Skew-2} yields
\begin{equation}\label{eq:L:Auto-2}
{\rm Im} \sum_{j \geq 1} \rho_{j}(AU_{0}\phi_{j}|LU_{0}\phi_{j})
+ {\rm Im}\sum_{j \geq 1}\rho_{j}(V[U_{0}]U_{0}\phi_{j}|LU_{0}\phi_{j})
= 0
\end{equation}
for all such operators $L$.

At this point, in a first step, assume that the integer $k$ is such that condition \eqref{eq:RhojDist} is fulfilled.
Choosing $M$ and $L$ as in \eqref{eq:ChoiceML}, and proceding exactly as in the proof of Lemma~\lemref{lem:ADiag}, using the fact that $\rho_{k} - \rho_{n} \neq 0$, we conclude that 
$$
(Au_{k} + V[U_{0}]u_{k}|u_{n}) = 0,
$$
that is 
$$
Au_{k} + V[U_{0}]u_{k} \in {\rm span}\{u_{n} \; ; \; n\neq k\}^\perp = {\rm span}\{u_{k}\}.
$$ 
This means that 
$$ Au_{k} + V[U_{0}]u_{k} = \lambda_{k}u_{k} $$ 
for some $\lambda_{k} \in {\Bbb C}$, but since $A + V[U_{0}]$ is a  self-adjoint operator, as a matter of fact we have $\lambda_{k} \in {\Bbb R}$. 
\smallskip

Next assume that the integer $k$ is such that the coefficient $\rho_{k}$ has multiplicity $m \geq 2$, that is condition \eqref{eq:RhojMult} is satisfied.
Arguing as above, we choose the operators $M$ and $L$ as in \eqref{eq:ChoiceML-2}, and conclude that 
$$(Au_{k+\ell} + V[U_{0}]u_{k+\ell} |u_{n}) = 0$$
 for all $n \not\in \left\{k+j \; ; \; 0 \leq j \leq m-1 \right\}$, that is:
$$
Au_{k+\ell} + V[U_{0}]u_{k+\ell} \in \left({\rm span}\left\{ u_{n} \; ; \; n \neq k+j, \; 0 \leq j \leq m-1 \right\}\right)^\perp.
$$
This means that if we set $H_{k} := {\rm span}\left\{u_{k+i} \; ; \; 0 \leq i \leq m-1 \right\}$, and 
$$
A_{0} u := -\Delta u + \left({\widetilde V}_{0} + V[U_{0}] \right) u
$$ 
then $A_{0} : H_{k} \dans H_{k}$ is a self-adjoint operator on the finite dimensional space $H_{k}$. Therefore there exists a unitary operator $U_{k}$, acting on this space, such that if for $0 \leq \ell \leq m-1$ we set ${\widehat u}_{k+\ell} = U_{k}u_{k+\ell} = U_{k}U_{0}\phi_{k+\ell}$, we have $A_{0}{\widehat u}_{k+\ell}=\lambda_{k+\ell}{\widehat u}_{k+\ell}$ for some $\lambda_{k+\ell} \in {\Bbb R}$.

However, since $U_{k}$ is a unitary operator on $H_{k}$ we have
$$
\sum_{\ell = 0}^{m-1}|U_{0}\phi_{k+\ell}|^2 = 
\sum_{\ell = 0}^{m-1}|U_{k}U_{0}\phi_{k+\ell}|^2 
$$
and thus
$$
\sum_{j=k}^{k+m-1}\rho_{j}|U_{0}\phi_{j}|^2 = 
\rho_{k}\sum_{\ell = 0}^{m-1}|U_{0}\phi_{k+\ell}|^2 = 
\rho_{k}\sum_{\ell = 0}^{m-1}|U_{k}U_{0}\phi_{k+\ell}|^2. 
$$
This means that if we set ${\widetilde U}_{k}\phi_{j} := U_{0}\phi_{j} $ if $j \not\in \{k+\ell \; ; \; 0 \leq \ell \leq m-1$ and ${\widetilde U}_{j} = U_{k}U_{0}\phi_{j}$ if $k \leq j \leq k+m-1$, we have $V[U_{0}] = V[{\widetilde U}_{k}]$, and finally this implies that
$$
A_{0}{\widehat u}_{k+\ell} = -\Delta {\widehat u}_{k+\ell} + \left({\widetilde V}_{0} + V[{\widetilde U}_{k}]\right){\widehat u}_{k+\ell} = \lambda_{k+\ell}{\widehat u}_{k+\ell}.
$$
\endproof
\medskip

As we may see from the above analysis, when all the $\rho_{j}$'s are distinct, then any unitary operator $U_{0}$ which minimizes $J$, yields a solution to the Schr\"odinger--Poisson system. However in the general case, when some of the coefficients $\rho_{k}$ have multiplicity $m_{k} \geq 2$, it is possible that one has to impose a unitary transformation $U_{k}$ in the space 
$$
H_{k}:={\rm span}\{U_{0}\phi_{k+\ell} \; ; \; 0 \leq \ell \leq m_{k} - 1\}
$$
in order to obtain a solution (note that these unitary transformations do not change the value of $V[U_{0}]$).
 
In other words, one may find a unitary operator $U_{k}$ on $H_{k}$ such that if $A_{k} := A_{|H_{k}}$ is the trace of $A$ on $H_{k}$, the operator  $U_{k}^*A_{k}U_{k}$ is diagonal. Thus since $\rho_{k}=\rho_{k+\ell}$ for $0 \leq \ell \leq m-1$, if we denote by ${\widehat U}$ the unitary operator obtained through the composition of all such operators $U_{k}$ and $U_{0}$, one has $J({\widehat U}) = J(U_{0})$. More precisely, we can we state the following corollary, which ends the proof of Theorem \lemref{lem:SPDiag}:

\begin{corollary}\label{lem:SolSP-2}
Under the assumptions of Theorem \lemref{lem:SolSP}, let $U_{k}$ be given by Lemma~\lemref{lem:SPDiag} when $k \geq1$ is such that \eqref{eq:RhojMult} is satisfied. Define the operator ${\widehat U}_{0}$  by ${\widehat U}_{0}\phi_{k} = U_{0}\phi_{k}$ when $k$ satisfies \eqref{eq:RhojDist}, and 
$$
{\widehat U}_{0}\phi_{k+\ell} := U_{k}U_{0}\phi_{k+\ell}, \quad\mbox{for }\, 0 \leq \ell \leq m-1, \quad\mbox{when \eqref{eq:RhojMult} is satisfied}.
$$
Then ${\widehat U}_{0}$ belongs to ${\Bbb S}$, while $J({\widehat U}_{0}) = J(U_{0})$ and $V[{\widehat U}_{0}] = V[{\widehat U}_{0}]$. Moreover setting $u_{j} := {\widehat U}_{0}\phi_{j}$ for $j \geq 1$,  there exists $\lambda_{j} \in {\Bbb R}$ such that 
$$
-\Delta u_{j} + \left({\widetilde V}_{0} + V[{\widehat U}_{0}]\right) u_{j} = \lambda_{j} u_{j}, \quad u_{j}\in {\Bbb H}_{1}, \quad (u_{j}|u_{k}) = \delta_{jk}
$$
and moreover $V := V[{\widehat U}_{0}]$ satisfies
$$
-\Delta V = \sum_{j \geq 1} \rho_{j} |u_{j}|^2, \quad V \in H^1_{0}(\Omega).
$$
\end{corollary}


\section{Further remarks}\label{sec:Remarks}

The one dimensional case $d=1$ is particularly simple to handle, using a completely different method. Indeed we point out that for a given potential $V \in C([0,1])$ the spectral sequence $(\lambda_{k},\phi_{k})_{k\geq 1}$ 
$$-\phi''_{k}  + (V + {\widetilde V})\phi_{k} = \lambda_{k}\phi_{k}, \quad \phi_{k}(0) = \phi_{k}(1) = 0, \quad \phi'_{k}(0) >0,$$
is well defined, and $\int_{0}^1\phi_{k}\phi_{j}dx = \delta_{kj}$, each eigenvalue $\lambda_{k}$ being simple. Using the simplicity of the eigenvalues, it is known that the mapping $V \mapsto \phi_{k}$ is continuous from $C([0,1])$ into $L^2(0,1)$ (see J. P\"oschel \& E. Trubowitz~\cite{PoTr1897}). If the coefficients $(\rho_{j})_{j \geq1}$ satisfy
$$\rho_{j} >0, \qquad \sum_{j \geq 1}\rho_{j} =: M < \infty$$
one can easily see that the mapping $F : C([0,1]) \dans L^1(0,1)$
$$V \mapsto F(V) := \sum_{k \geq 1}\rho_{k}|\phi_{k}|^2$$
is continuous.

Now consider the mapping $B : L^1(0,1) \dans L^1(0,1)$ defined by $Bf := v$ where $v\in C^1([0,1])$ is given by
$$-v'' = f \quad \mbox{in }\, (0,1), \qquad v(0) = v(1) = 0.$$
Clearly $B$ can also be considered as a linear mapping on $C([0,1])$, and $B$ is compact. Also observe that the second equation in \eqref{eq:2} is equivalent to find $V \in C([0,1])$ such that
$$V - B F(V) = 0.$$
Denoting by $T(V) := BF(V)$, we know that $I - T$ is a compact perturbation of the identity on $C([0,1])$ and one can check easily that there exists $R>0$ such that for any $\theta \in [0,1]$
$$V - \theta T(V) = 0 \imply \|V\|_{\infty} < R.$$
Thus the invariance by homotopy of the Leray--Schauder topological degree implies that for all $\theta \in [0,1]$ we have
$${\rm deg}(I-\theta T, B(0,R),0) = {\rm deg}(I, B(0,R),0) = 1,$$
which means in particular that ${\rm deg}(I- T, B(0,R), 0) = 1$. Therefore there exists at least one $V \in C([0,1])$ such that $V - T(V) = 0$, that is the system \eqref{eq:1}--\eqref{eq:2} has at least one solution, when $\Omega = (0,1)$.

We point out also that in the case in which $\rho_{m}$ is a function of $\lambda_{m}$, for instance $\rho_{m} := \exp(-\lambda_{m})$ (see F.~Nier \cite{Nier1993CPDE}), the same approach can be applied.
In dimensions $N \geq 1$, J.Ph.~Solovej \cite{Solovej:1991} considers a {\it one-particle density matrix\/} defined by $f\mapsto \gamma f := \sum_{m\geq 1}\rho_{m}(u_{m}|f)u_{m}$ and minimizes a functional depending on $\gamma$ under the constraint ${\rm spectrum}(\gamma) = \{\rho_{m} \; ; m \geq 1\}$.
See also E.~Prodan \cite{Prodan:2005}, E.~Prodan \& P.~Nordlander \cite{ProdanNordlander:2001} where Hartree-Fock approximations are considered.



\bibliographystyle{plain}

\end{document}